\documentclass[a4paper,10pt]{article}
\usepackage[utf8x]{inputenc}
\usepackage{amsmath,amsfonts,theorem}
\title{New example of modified moduli space of  special Bohr - Sommerfeld lagrangian submanifolds.}
\author{Nik. A. Tyurin\footnote{ The author is partially supported by Laboratory of Mirror Symmetry NRU HSE, RF Government grant, ag. 14.641.31.0001 }}

\begin{document}

\maketitle

\begin{center}
 {\em  Bogolyubov Lab. of Theoretical Physics, JINR (Dubna), \\  Lab. of Mirror Symmetry and Automorphic forms, HSE (Moscow)}
 \end{center}

\begin{abstract} We present an example of modified moduli space of special Bohr - Sommerfeld lagrangian submanifolds for the case
when the given algebraic variety is the full flag $F^3$ for $\mathbb{C}^3$ and the very ample bundle is $K^{- \frac{1}{2}}_{F^3}$.

\end{abstract}

Let $X$ be a compact simply connected algebraic variety of complex dimension $n$, $L \to X$ be a very ample line bundle and $h$ be an appropriate hermitian structure on $L$ such that
it gives a Kahler structure on $X$ with the Kahler form $\omega_h$. Let $\vert L \vert$ be the projective space corresponding to holomorphic section space for $L$.
Fix a topological type ${\rm top} S$ and a middle class $[S] \in H_n (X, \mathbb{Z})$ where $S$ is a real $n$ - dimensional orientable manifold.

Consider the space of pairs $(\{ S \}, D)$ where 

--- $D \in \vert L \vert$ is a zero divisor of a holomorphic section from $H^0(X, L)$;

--- $\{ S \}$ is a class of smooth  homologically non trivial D- exact with respect to $D$ lagrangian submanifolds from the complement $X \backslash D$ modulo Hamiltonian isotopies in $X \backslash D$,
which are of topological type ${\rm top} S$ and represent the class $[S]$ being considered in $X$.

The definition of D - exact lagrangian submanifolds in full generality is the following. 

{\bf Definition.} {\it Let $(X, \omega, D)$ is a simply connected symplectic manifold with integer symplectic form $\omega$ and $D \subset X$ be 
 $2n-2$ submanifold whose homology class is Poincare dual to $[\omega]$. Then we say that a lagrangian submanifold $S \subset X$ is D - exact iff
for any loop $\gamma \subset S$ and any disc $B_2 \subset X$ bounded by $\gamma$ the topological intersection $B_2 \cap D$ equals to $\int_{B_2} \omega$.}

Note that this property  does not depend on any additional structures --- just on the mutual arrangement of $S$ and $D$. 
On the other hand this defintion can be extended even to the case when $S \cap D \neq \emptyset$. It is not hard to see that this definition is applicable in the case we presented above.

Coming back to the definition of the modified moduli space.
    Step by step, first we take a divisor $D \in \vert L \vert$ and consider the complement $X \backslash D$. In this complement we have the group $H_n(X \backslash D, \mathbb{Z})$ with the natural
epimorphism $\pi: H_n(X \backslash D, \mathbb{Z}) \to H_n(X, \mathbb{Z})$. Take the preimage $\pi^{-1}([S])$ and for each class from $\pi^{-1}([S]) = \{\chi_1, ... \chi_m\},
\chi_i \in H_n(X \backslash D, \mathbb{Z})$, realizable by smooth D - exact lagrangian submanifolds,  find  such smooth D- exact with respect to $D$ lagrangian submanifolds which 
represent in the complement $X \backslash D$ the classes
$\chi_i$. These lagrangian submanifolds form  spaces of solutions ${\cal L}_{D-ex}^{\chi_i}$ for each $i$. The last (but not least) step is the factorization
of each ${\cal L}^{\chi_i}_{D-ex}$ by the natural action of Hamiltonian isotopies on the complement $X \backslash D$. Each element of the factorized space gives a point
$(\{S \}, D)$ of our modified moduli space $\tilde {\cal M}_{SBS}$. The motivation and some prehistory can be found in [1].

Let $X = F^3$ be the full flag variety in $\mathbb{C}^3$, realized as a hypersurface $X = \{ \sum_{i=0}^2 x_i y_i = 0 \}$ in the direct product $\mathbb{C} \mathbb{P}^2_x \times \mathbb{C} \mathbb{P}^2_y$
where $[x_i], [y_j]$ are homogenious coordinates on each direct summand. The line bundle ${\cal O}(1,1)$ being restricted to $X$ is very ample, and the corresponding symplectic form
$\omega$ is given by the direct sum of the lifted standard Fubini - Study forms from the direct summands. Thus we consider $L = {\cal O}(1,1)|_X$.

It is known that the flag variety $F^3 = X$ with the standard symplectic form $\omega$ contains a lagrangian 3 - sphere which is called {\it the Gelfand - Zeytlin} sphere since it arose in the framework of the Gelfand - Zeytlin systems. It is explicitly presented in the homogenious coordinates by the condition
$$
S_{GZ} = \{ [x_0: x_1: x_2] \times [y_0 = \bar x_0: y_1 = \bar x_1: y_2 = - \bar x_2] \vert \sum_{i=0}^2 x_i y_i = 0 \} \subset X.
\eqno (*)
$$
It is not hard to check that the condition above describes a lagrangian sphere: since the antidiagonal embedding of $\mathbb{C} \mathbb{P}^2$ to the direct product
$\mathbb{C} \mathbb{P}^2_x \times \mathbb{C} \mathbb{P}^2_y$ is lagrangian, the subset $\{  [x_0: x_1: x_2] \times [y_0 = \bar x_0: y_1 = \bar x_1: y_2 = - \bar x_2] \}$
is Hamiltonian isotopic to the antidiagonal embedding thus is lagrangian as well, and the intersection of this subset and $X$ is {\it cotransversal} one deduces
that it is lagrangian. On the other hand $S_{GZ}$ is described by the condition $\vert x_0 \vert^2 + \vert x_1 \vert^2 - \vert x_2 \vert^2 = 0$ in $\mathbb{C}
\mathbb{P}^2$ therefore it is isomorphic to 3 - sphere. 

Since the same is true for the case when minus in the formula (*) is placed not before $\bar x_0$ but before any other $\bar x_i$ we get
another lagrangian sphere of the same type, namely
$$
S_0 = \{[x_0: x_1: x_2] \times [- \bar x_0: \bar x_1: \bar x_2]\}, S_1 = [x_0: x_1: x_2] \times [\bar x_0: - \bar x_1: \bar x_2] \},
$$
and $S_2 = S_{GZ}$ --- all are lagrangian spheres, which are Hamiltonian isotopic in $X$. 
At the same time the homology class $[S_i] \in H_3 (X, \mathbb{Z})$ is trivial since $H_3 (X, \mathbb{Z})$ is trivial itself.

We would like to describe the moduli space $\tilde {\cal M}_{SBS}$ for the following topological data: our $S$ is isomorphic to 3 - sphere homologically trivial in $X$. Note that 
in this case {\it every} lagrangian sphere must be exact having trivial fundamental group, which drastically simplify the analysis. 

We start with an irreducible divisor  $D \in \vert L \vert$ and our aim is to find the classes of lagrangian spheres in the complement $X \backslash D$ up to Hamiltonain isotopies.
For this first study the homotopy type of the complement $X \backslash D$. Being a hypersurface in the direct product $\mathbb{C} \mathbb{P}^2_x \times \mathbb{C} \mathbb{P}^2_y$ our
flag variety admits the canonical projection to the first summand $\pi_x : X \to \mathbb{C} \mathbb{P}^2_x$, whose fibers are projective lines in the second projective plane.
Any irreducible $D \subset X$ gives a section of this projection plus three fibers, therefore the complement $X \backslash D$ is given by a complex line bundle over
$\mathbb{C} \mathbb{P}^2_x \backslash \{p_1, p_2, p_3 \}$ where $p_i$ is a point in this projective plane. To illustrate it we can take say the following $D \subset X$:
$$
D = \{x_0 y_0 - x_1 y_1 + i x_2 y_2 = 0, x_0 y_0 + x_1 y_1 + x_2 y_2 =0  \} \subset X.
$$
It is equivalent to conditions $(1-i) x_0 y_0 = (1+i) x_1 y_1  = - x_2 y_2$ therefore for each $[x_0: x_1: x_2]$ we either have exactly one point in $\pi^{-1}_x([x_0: x_1: x_2])$
if at least two coodrinates are non zero or whole the fiber if two coordinates are zero. Therefore $X \backslash D$  is naturally isomorphic to
the product 
$$
(\mathbb{C} \mathbb{P}^2_x \backslash \{[1:0:0], [0:1:0], [0: 0: 1] \}) \times \mathbb{C}.
$$
Contracting the punctured fibers we get that $X \backslash D$ is homotopic to the complement  $\mathbb{C} \mathbb{P}^2_x \backslash \{[1:0:0], [0:1:0], [0: 0: 1] \}$.
Note that for any irreducible divisor the picture is essentially the same. 

Exclusion of three distinct points from $\mathbb{C} \mathbb{P}^2_x$ generates non trivial homotopy group $\pi_3(\mathbb{C} \mathbb{P}^2_x \backslash \{p_1, p_2, p_3 \})$
equals to $\mathbb{Z} \oplus \mathbb{Z}$; thus $H_3 (X \backslash D, \mathbb{Z})$ is the same $\mathbb{Z} \oplus \mathbb{Z}$, and one can expect that certain classes
from this lattice can be realized by smooth lagrangian spheres. 

This is indeed the case: we claim that for each $i = 0, 1, 2$ lagrangian sphere $S_i$ defined above:

i) does not intersect divisor $D$;

ii) lies in a non trivial class $a_i \in H_3(X \backslash D, \mathbb{Z})$ such that $a_i \neq a_j$ for $i \neq j$;

iii) no other classes of lagrangian spheres do exist.

First check item i): since our divisor is given by the conditions $(1-i)x_0 y_0 = (1+i) x_1 y_1 = - x_2 y_2$ substitute there $y_0 = - \bar x_0, y_1 = \bar x_1, y_2 = \bar x_2$ as it is
for $S_0$. This leads to the condition $(1-i)\vert x_0 \vert^2 = (1+i) \vert x_1 \vert^2 = - \vert x_2 \vert^2$ which can take place if and only if $x_0 = x_1 = x_2 =0$ but such a point does not exist
on the projective plane. Therefore $D \cap S_0 = \emptyset$. The same arguments work for $S_1$ and $S_2$ as well. 

Further, consider the projections of $S_i$ under $\pi_x$ to the first projective plane. Note that the projection realizes the contraction of $X$ to $\mathbb{C} \mathbb{P}^2_x$;
the image of $S_0$ is given by the equation $\pi_x(S_0) = \{- \vert x_0 \vert^2 + \vert x_1 \vert^2 + \vert x_2 \vert = 0 \} \subset \mathbb{C} \mathbb{P}^2_x$, it is a smooth 3 - dimensional sphere
which divides the projective plane into two open parts. These parts are described by the sign of the value of the real function
$$
F_0 = \frac{- \vert x_0 \vert^2 + \vert x_1 \vert^2 + \vert x_2 \vert^2}{\sum_{i=0}^2 \vert x_i \vert^2},
$$
which maps $\mathbb{C} \mathbb{P}^2_x$ to the segment $[-1, 1]$. Thus function has non degenerated critical point at $[1:0:0]$ and degenerated critical set at the line $x_0 = 0$.
The value $F_0 = 0 $ is non critical, therefore $S_0$ is a smooth sphere. At the same time at points $[1:0:0]$ and $[0:1:0], [0:0:1]$ the function $F_0$ has different signs,
consequently $S_0$ lies in a non trivial class in $\pi_3(X \backslash D)$ since $\pi_x$ realizes a homotopy $X \backslash D \mapsto \mathbb{C} \mathbb{P}^2_x \backslash \{[1:0:0], [0:1:0],
[0:0:1]\}$ and $\pi_x(S_0)$ represents a non trivial class. 

The same arguments show that $S_1$ does represent a non trivial class as well, but it is a different class from $H_3(X \backslash D, \mathbb{Z})$ since
it is projected to a different class for the punctured projective plane. Indeed, $S_0$ is ``centered'' in $[1:0:0]$ and the resting points lie ``outside'' of $S_0$;
at the same time $S_1$ is ``centered'' in $[0:1:0]$ and the resting points $[1:0:0], [0:0:1]$ lie ``outside'' of it; the same is true for $S_2$. The basis in
$H_3(X \backslash D, \mathbb{Z})$ can be choosen in such a way that $S_0, S_1, S_2$ represent the classes $(1,0), (0,1)$ and $(1,1)$ respectively.     
In particular this fact ensures us that $S_i$ and $S_j$ are not Hamiltonian isotopical for $i \neq j$.

For momentary prove of iii) we need that $S_{GZ}$ is essentially unique lagrangian sphere in $F^3$ up to Hamiltonian isotopy
since it were possible to establish the fact without any references.

Thus we can say that for irreducible divisors the classes can be described by triple of points $p_1, p_2,  p_3 \in \mathbb{C} \mathbb{P}^2_x$; below we present an algebraic way to
define this attachment. Now we consider reducible divisors from the geometrical poitn of view. 

A reducible divisor $D \in \vert L \vert$ is presented by pair of sections of ${\cal O}(1)$ on each copy of the projective planes: therefore we can represent it by two projective lines
$l_x \subset \mathbb{C} \mathbb{P}^2_x, l_y \subset \mathbb{C} \mathbb{P}^2_y$, and the divisor $D \subset X$ is given by the union of $\pi_x^{-1}(l_x) \cup \pi_y^{-1} (l_y)$.
These two components of $D$ in $X$ are isomorphic to the Hirzebruch surface $F_1$, and if we study the projection $\pi_x: X \to \mathbb{C} \mathbb{P}^2_x$ then
the first component is given by the fibers over $l_x$ while the second component presents a section of this projection plus the fiber over a point which corresponds to
the porjective line $l_y$ (thus it is isomorphic to the projective plane with one point blown up). Therefore we have two different situations:
in more general case the point, corresponding to $l_y$ does not lie on $l_x$ in $\mathbb{C} \mathbb{P}^2_x$, in more specific case the point lies on $l_x$.

In the first case the homotopy type of $X \backslash D$ is the same as of $\mathbb{C}^2 \backslash \{pt \}$: indeed, after the cancelation of the first component
we get a projective bundle over $\mathbb{C} \mathbb{P}^2_x \backslash l_x$ which means the projective bundle over $\mathbb{C}^2$; then the cancelation of the second component
removes one point in each fiber and removes totally one fiber over the point which corresponds to $l_y$. Therefore $X \backslash D$ is isomorphic
to $(\mathbb{C}^2 \backslash \{ pt \}) \times \mathbb{C}$, and consequently $H_3(X \backslash D, \mathbb{Z}) = \mathbb{Z}$. Now we can again show that
there exists unique up to Hamitonian isotopy smooth lagrangian sphere in the complement $X \backslash D$, which is projected by $\pi_x$ to
a 3 - sphere in $\mathbb{C}^2 \backslash \{ pt \}$ ``centered'' exactly in this point. This point $p$ characterizes the sphere, and we attach to this point
the corresponding class of lagrangian spheres up to Hamiltonian isotopy.

The last case when $l_x$ and $l_y$ are related by the condition that the point in $\mathbb{C} \mathbb{P}^2_x$, corresponding to $l_y$, lies on $l_x$,
gives trivial set of lagrangian spheres since in this case $X \backslash D$ is isomorphic to $\mathbb{C}^2 \times \mathbb{C}$ therefore $H_3(X \backslash D, \mathbb{Z})$
is trivial. For this case every Gelfand - Zeytlin sphere intersects the divisor.

Summing up, we see that the moduli space $\tilde {\cal M}_{SBS}$ can be described as a subset of the direct product $\vert L \vert \times \mathbb{C} \mathbb{P}^2_x$:
for each element of $\vert L \vert \cong \mathbb{C} \mathbb{P}^7$ the classes of smooth lagrangian spheres are given by the ``centers'' of
their projections to $\mathbb{C} \mathbb{P}^2_x$; and the complete linear system $\vert L \vert$ is stratified by the conditions ``three points in $\mathbb{C} \mathbb{P}^2_x$'',
``one point in $\mathbb{C} \mathbb{P}^2_x$'' and `` no points in $\mathbb{C} \mathbb{P}^2_x$. Note that the picture looks quite similiar to the 
example from [1], for $X = \mathbb{C}\mathbb{P}^1, L = {\cal O}(3)$.

The detailed analysis shows that the answer even closer to the result for $X = \mathbb{C}\mathbb{P}^1, L = {\cal O}(3)$. 

Study the situation algebracally: with respect to the fixed homogenious coordinates $[x_0: x_1: x_2], [y_0: y_1: y_2]$ every divisor $D$ from $\vert L \vert$
is given by two equations: $\sum_{i=0}^2 x_i y_i = 0, \sum_{i,j} a_{ij} x_i y_j = 0$ where in the last expression $i$ and $j$ are taken from 0,1 or 2.
The representation of $D$ by the numbers $\{ a_{ij} \}$ is not unique since we can add any $\lambda$ to $a_{00}, a_{11}$ and $a_{22}$ without changing
of the system and as well we can scale the numbers $a_{ij}$ therefore we normalize the matrix
$$
A = ( a_{ij} )
$$
such that ${\rm tr} A = 0$ and consider it up to scaling. 

For a given divisor $D$ write the defining system as
$$
\begin{aligned}{1}
  x_0 y_0 + x_1 y_1 + x_2 y_2 = 0 \\
(a_{00} x_0 + a_{10}x_1 + a_{20}x_2) y_0 + (a_{01} x_0 + a_{11} x_1 + a_{21} x_2) y_1 \\ + (a_{02} x_0 + a_{12} x_1 + a_{22} x_2) y_2 =0 \\
\end{aligned}
$$
For a particular choice of $[x_0: x_1: x_2]$ understood as parameters the last system looks like a system of linear equations in three variables $y_0, y_1, y_2$,
which has a unique up to scale non zero solution if and only if the system
$$
\begin{aligned}{1}
 (a_{00} x_0 + a_{10}x_1 + a_{20}x_2) = \lambda x_0\\
 (a_{01} x_0 + a_{11}x_1 + a_{21}x_2) = \lambda x_1\\
(a_{02} x_0 + a_{12}x_1 + a_{22}x_2) = \lambda x_2 \\
\end{aligned}
\eqno (**)
 $$
does not admit non zero solution for any $\lambda$. It implies that  any divisor $D$ represents a section of the projection $\pi_x: X \to \mathbb{C} \mathbb{P}^2_x$
over the subset $N_A \subset \mathbb{C} \mathbb{P}^2_x$ consists of such $[x_0: x_1: x_2]$ that the system (**) does not admit non zero solutions. But it is exactly
the condition that vector $(x_0, x_1, x_2)$ is not  an eigenvector for the matrix $A$ with an eigenvalue $\lambda$. It follows that the divisors from the complete linear system
$\vert L \vert$ can be combined with respect to the following stratified conditions:

1) (generic case) the matrix $A$ has three eigenvectors with  different eigenvalues $\lambda_1, \lambda_2, \lambda_3$, $\lambda_i \neq \lambda_j$, therefore we have 
three distinct points $p_1, p_2, p_3 \in \mathbb{C} \mathbb{P}^2_x$
--- and thus three different classes of smooth lagrangian spheres;

2) the matrix $A$ has one multiple eigenvalue, say $\lambda_1 = \lambda_2$, and admits one eigenvector and one eigenspace, therefore we have one distinct point $p$ and the corresponding line $l_x$
in the projective space $\mathbb{C} \mathbb{P}^2_x$ such that $p$ does not lie on $l_x$, --- and for such a divisor we have one class of lagrangian spheres;

3) the matrix $A$ has one multiple eigenvalue, say $\lambda_1 = \lambda_2$, but admits two eigenvectors and no eigenspace (the case of Jordan cell), therefore
after the projection by $\pi_x$ we get $\mathbb{C} \mathbb{P}^2_x \backslash \{p_1, p_2 \}$ but $\pi_3 (\mathbb{C} \mathbb{P}^2_x \backslash \{p_1, p_2 \}) = \mathbb{Z}$,
--- therefore we again as in the case 2) get only one class of lagrangian spheres;

4) the matrix $A$ has $\lambda_1 = \lambda_2 = \lambda_3 = 0$ and is presented by one Jordan cell, so it admits a single eigenvector, --- this corresponds
to  absence of lagrangian spheres  in $X \backslash D$ since $\pi_3(\mathbb{C} \mathbb{P}^2_x \backslash \{pt \}) = 0$;

5) the matrix $A$ has $\lambda_1 = \lambda_2 = \lambda_3 = 0$ and contains $2 \times 2$ Jordan cell, thus we have just an eigenspace of dimension 2 and no
separate eigenvector, so this is the case discussed above when $l_x$ and $l_y$ are related, and --- in this case no lagrangian spheres in the complement.

 Summing up the cases 1) - 5) we see that the classes of lagrangian spheres correspond to {\it single} eigenvalues of the matrix $A$. Indeed, in the case 1) we have 
three single eigenvalues --- and three classes of lagrangian spheres; in the cases 2) and 3) we have one single eigenvalue $\lambda_3$ --- and unique class of lagrangian sphere;
at finally in the cases 4) and 5) we do not have lagrangian spheres in the complement $X \backslash D$.

These arguments leads to the following realization of the modified moduli space $\tilde {\cal M}_{SBS}$: in the direct sum $H^0(X, {\cal O}(1,1)|_X) \oplus \mathbb{C} =
\mathbb{C}^8 \oplus \mathbb{C}$ consider an affine  hypersurface  given in the affine coordinates $(a_{ij}, z) \vert a_{00}+ a_{11} + a_{22} = 0$ by the cubic equation
${\rm det} (A - zI) = 0$ which is {\it homogenious} since the multiplication of matrix $A$ by a constant leads to the multiplication of the eigenvalues by the same constant.
Therefore this equation defines a projective cubic hypersurface $Y \subset \mathbb{P}(H^0(X, {\cal O}(1,1)|_X) \oplus \mathbb{C}) = \mathbb{C} \mathbb{P}^8$ definig a finite covering
of the complete linear system $\vert L \vert$. The ramification divisor $\Delta \subset Y$ is described by the following condition: since ${\rm det} (A - zI) = - z^3 - a z + b$
where $a$ has degree 2 in the homogenious coordinates in $\vert L \vert$ and $b = {\rm det} A$ has degree 3, then the multiple eigenvalues correspond to non trivial solutions of
the system
$$
\begin{aligned}{2}
 z^3 + az - b = & 0, \\
(z^3 + az - b)'_{z} = & 3 z^2 + a = 0. \\
\end{aligned}
 $$
Therefore $\Delta \subset Y \subset \mathbb{C} \mathbb{P}^8$ is a divisor from the complete linear system $\vert {\cal O}(2)|_Y \vert$ (note that $3z^2 + a = 0$ is homogenious
of degree 2 in   $H^0(X, {\cal O}(1,1)|_X) \oplus \mathbb{C}$).

{\bf Proposition.}  {\it The modified moduli space for the lagrangian 3- spheres in the flag variety $F^3$ and the line bundle $ K^{-\frac{1}{2}}_{F^3}$ is isomorphic to
$$
\tilde {\cal M}_{SBS} = Y \backslash \Delta.
$$
In particular it admits a natural compactification, isomorphic to $Y$}.

It is very interesting result since we get again the geometrical data of the same form (algebraic variety, very ample divisor). In particular we can
attach to our given pair $(F^3, K_{F^3}^{-\frac{1}{2}})$ either the Weinstein sceleton of $Y \backslash \Delta$ or thew finite set of the smooth exact lagrangian
submanifolds of $Y \backslash \Delta$ up to Hamiltionian isotopies.    

$$
$$

{\bf Bilbliography.}

[1] {\bf Nik. A. Tyurin}, {\it ''Towards the moduli space of special Bohr - Sommerfeld lagrangian cycles``}, {\bf arXiv:1708.00723}.

\end{document}